\documentclass{article}

\usepackage[english]{babel}
\usepackage{amsmath}
\usepackage{amsthm}
\usepackage{amssymb}

\newtheorem{thm}{Theorem}

\newtheorem{lem}{Lemma}

\newtheorem{prop}{Proposition}

\theoremstyle{definition}

\theoremstyle{remark}

\theoremstyle{remark}

\def\cM{\mathcal M}

\def\cP{\mathcal P}

\def\cH{\mathcal H}

\def\bR{\mathbb R} 

\def\bS{\mathbb S}

\def\bC{\mathbb C}

\def\NN{\mathbf N}

\def\Ker{\operatorname{Ker}}
\def\krat{\;\rfloor}

\def\La{\Lambda}

\def\pa{\partial}

\begin{document}

\title{The Fischer decomposition for Hodge-de Rham systems in Euclidean spaces}

\author{Richard Delanghe, Roman L\' avi\v cka and Vladim\'ir Sou\v cek}

\maketitle

\begin{abstract}
The classical Fischer decomposition of spinor-valued
polynomials is a~key result on solutions of the Dirac
equation in the Euclidean space $\bR^m.$ As is well-known, it can be understood as an irreducible decomposition with
respect to the so-called $L$-action of the Pin group $Pin(m).$ 
But, on Clifford algebra valued polynomials, we can consider also the $H$-action of $Pin(m).$ 
In this paper, the corresponding Fischer decomposition for the $H$-action is obtained.
It turns out that, in this case, basic building blocks are the spaces of homogeneous solutions to the Hodge-de Rham system.
Moreover, it is shown that
the Fischer decomposition for the $H$-action can be viewed even as a refinement of
the classical one.

\medskip\noindent{\bf Keywords:} Fischer decomposition, Clifford analysis, Hodge-de Rham equation, spherical monogenics

\medskip\noindent{\bf Mathematics Subject Classification (2000)} 30G35, 58A10
\end{abstract}


\section{Introduction}
Clifford analysis is, from the very beginning, considered to be a refinement of harmonic analysis
for Clifford algebra (or spinor) valued functions. The perfect description of this statement is
the Fischer decomposition of spinor-valued polynomials.

Let us first recall the Fischer decomposition of the space
$\cP$ of complex-valued polynomials in the Euclidean space $\bR^m.$ 
Denote by $\cH_k$ the space of $k$-homogeneous harmonic polynomials in $\bR^m.$ 
Then, under a~natural action of the orthogonal group $O(m),$ the space 
$\cP$ has an irreducible (not multiplicity free) decomposition
\begin{equation}\label{Fischer}
\cP=\bigoplus_{k=0}^{\infty}\bigoplus_{p=0}^{\infty}r^{2p}\cH_k
\end{equation}
where $r^2=x_1^2+\cdots +x_m^2$ for the vector variable $\underline{x}=(x_1,\ldots,x_m)\in\bR^m.$

For spinor-valued polynomials,
there is a refinement of this decomposition.
Let $\bC_m$ be the complex Clifford algebra generated by vectors of the standard basis $(e_1,\ldots,e_m)$ of $\bR^m.$
Recall that the $L$-action of the Pin group $Pin(m)$ on functions $f:\bR^m\to\bC_m$ is defined by
\begin{equation}\label{Laction}
[L(s)(f)](\underline{x}) = s\,f(s^{-1}\underline{x}s),\ s\in Pin(m)\text{\ \ and\ \ }\underline{x}=(x_1,\ldots,x_m)\in\bR^m.
\end{equation}
Denote by $\bS$ a~basic spinor representation for $Pin(m).$ 
As is well-known, the spinor space $\bS$ can be realized inside the Clifford algebra $\bC_m.$ 
Let us remark that, under the $L$-action, the space $\cP(\bS)=\cP\otimes\bS$ of spinor valued polynomials forms a~$Pin(m)$-module.
Denote by $\cM_k(\bS)$ the space of $k$-homogeneous polynomials $P\in\cP(\bS)$ which are (left) monogenic,  
i.e., which satisfy the Dirac equation $\underline{D}P=0$ where $$\underline{D}=e_1\pa_{x_1}+\cdots+e_m\pa_{x_m}.$$ 
Now we are ready to state the Fischer decomposition (sometimes called also Almansi decomposition) for this case. Namely, under the $L$-action, the space 
$\cP(\bS)$ has an irreducible (not multiplicity free) decomposition 
\begin{equation}\label{FischerPin}
\cP(\bS)=\bigoplus_{k=0}^{\infty}\bigoplus_{p=0}^{\infty}\underline{x}^p\cM_k(\bS)
\end{equation}
with $\underline{x}=e_1x_1+\cdots+e_mx_m.$ See \cite{BSES,MR,rya}.
As $\cH_k\otimes\bS=\cM_k(\bS)\oplus\underline{x}\cM_{k-1}(\bS)$ and $\underline{x}^2=-r^2,$
it is easy to see that (\ref{FischerPin}) is a~real refinement of (\ref{Fischer}).

The main aim of the underlying paper is to show that
there exists a natural further refinement of the monogenic Fischer decomposition \eqref{FischerPin}. It is 
quite surprising that such a finer Fischer decomposition was not described earlier.
It was the study of special solutions of the Dirac equation which led to such a~refinement.
By special solutions we mean just solutions having their values
in a chosen subspace $V$ of the Clifford algebra $\bC_m$. There are a lot of possibilities for a choice 
of $V,$ but it is clearly preferable to choose the subspace $V$ having some special properties. 

Typical examples
are solutions of the Dirac equation having values in  spinor
subspaces of the Clifford algebra. 
This case is closely related to the 
$L$-action \eqref{Laction}.
Indeed, it is well-known that the Clifford algebra $\bC_m$, considered as a 
$Pin(m)$-module by left multiplication, decomposes into many equivalent spinor submodules. 
Moreover, for every choice of the spinor
submodule, the $Pin(m)$-module of spinor-valued solutions has quite analogous properties.
 
Another interesting example of special solutions of the Dirac equation is given by the so-called generalized Moisil-Th\'eodoresco system (GMT system for short).
A~lot of interest has recently been paid to GMT systems (see \cite{DLS} and the references there).
In this case, the space $V$ is supposed to be invariant under another (both side) action of the Pin group, namely the so-called $H$-action.
The $H$-action on Clifford algebra valued functions $f:\bR^m\to\bC_m$ is given by
\begin{equation}\label{Haction}
[H(s)(f)](\underline{x}) = s\,f(s^{-1}\underline{x}s)s^{-1},\ s\in Pin(m)\text{\ \ and\ \ }\underline{x}\in\bR^m.
\end{equation}
In what follows, we shall use the language of differential forms.  
Indeed, following \cite{BDS}, we identify naturally the Clifford algebra $\bC_m$ with the Grassmann algebra $\La^*(\bC^m)$ 
and we study the space $\cP^*=\cP\otimes\La^*(\bC^m)$ of polynomial differential forms instead of Clifford algebra valued polynomials.
Then the $H$-action translates into a~natural action of the orthogonal group $O(m)$ on $\cP^*$ and the Dirac operator $\underline{D}$ corresponds to the operator $d+d^*.$ Here $d$ and $d^*$ are, respectively, the standard de Rham differential and its adjoint (see \eqref{dd^*} below). 
As an $O(m)$-module, the 
Grassmann algebra $\La^*(\bC^m)$ has a multiplicity free irreducible decomposition
$$\Lambda^*(\bC^m)=\bigoplus_{s=0}^m\Lambda^s(\bC^m)$$
with $\Lambda^s(\bC^m)$ being the space of $s$-vectors over $\bC^m.$ A~GMT system is then defined as the homogeneous system obtained by restricting the operator $d+d^*$ to functions having values in the space
$$V=\bigoplus_{s\in S}\Lambda^s(\bC^m)$$
for some (suitable) subset $S\subset\{0,1,\ldots,m\},$ i.e.
$$(d+d^*)P=0\mbox{\ \ \ for $V$-valued $P.$}$$
In particular, for $V=\Lambda^s(\bC^m),$ the corresponding GMT system coincide with
the so-called Hodge-de Rham system
\begin{equation}\label{Hodge}
d P=0,\, d^* P=0.
\end{equation}

Various versions of GMT systems of PDE's were studied for a long time (in particular in low dimensions) and they 
were used in many different applications. 
Applications in numerical analysis and engineering sciences can be found in \cite{gs1,gs2}.
In a~review paper \cite{spr}, you can find various generalizations of the well-known Hodge-de Rham decomposition of smooth 1-forms, including decompositions for quaternionic and Clifford algebra valued functions and the Almansi (i.e., monogenic Fischer) decomposition. For applications in theory of electromagnetic fields we can refer to \cite{S_Maxwell}.

The Fischer decomposition always played a key role in Clifford analysis. 
In \cite{lav_Fischer}, the Fischer decomposition for the $H$-action has been recently applied to inframonogenic functions introduced in \cite{MPS}. Moreover, in \cite{lav_Fischer}, the obtained results for the space $\cP^*$ are translated back into the framework of Clifford analysis. 
For yet another application, we can refer to \cite{lavSL3}.

Recently, the Fischer decomposition (together with the Cauchy-Kovalevskaya extension) was used systematically for construction of orthogonal bases in the spaces of homogeneous polynomial solutions.
In the classical Clifford analysis, it has a quite long history (see \cite{BGLS} for historical account,
various results can be found in \cite{Bock2010c,Bock2010a,Bock2009,BG,CacGueMal,CacGueBock,BockCacGue,cac,CM06,CM07,CM08,
FCM,FM,Gurlebeck1999,NGue2009,lavSL2,step2,Malonek1987,MS,mor09,som,van}).
Analogous results in Hermitean Clifford analysis are described in \cite{ckH, kerH,GTinH, GT2H}.
Finally, in \cite{DLS4}, the Fischer decomposition for the $H$-action plays a~key role in constructing orthogonal bases not only for the spaces $H^s_k$ of solutions to the Hodge-de Rham system but even for the spaces of homogeneous solutions of an arbitrary generalized
Moisil-Th\'eodoresco system.

In this paper, we will establish (using results from \cite{hom}) a form of the Fischer decomposition  appropriate for the $H$-action. The theory of the Howe duality developed in \cite{how} shows us that we may
expect a further refinement of the monogenic Fischer decomposition (\ref{FischerPin}), see \cite{DLS3} for details.
Indeed, this is visible from the form of invariants contained in the polynomial spaces considered.
For scalar valued functions, invariant polynomials are generated by powers of $r^2,$
and the basic equation is the Laplace equation.
For spinor-valued polynomials with the left action, we have to look for invariants with
values in the space of endomorphisms of the spinor space (which is, basically, the corresponding
Clifford algebra) and there is just a new invariant $\underline{x},$ acting as a refinement of $r^2.$
The corresponding basic equation is the Dirac equation.

In the case of the both side action, 
we deal with the space $\cP^*$ of $\La^*(\bC^m)$-valued
polynomials. The space of invariants with values in the space of endomorphisms of the Grassmann algebra $\La^*(\bC^m)$ is now much richer.  It is generated by two elements ${x}$ and ${x}^*$ which  correspond to the differential operators $d^*$ and $d$ by the Fischer duality (see \eqref{xx^*} below). 
Consequently, the corresponding basic system of equations is the Hodge-de Rham system (\ref{Hodge})
and the space of invariants consists of polynomials in $x$ and $x^*.$ Actually,  due to the fact
that ${x}^2=({x}^*)^2=0,$ such invariants are generated by the set 
\begin{equation}\label{Omega}
\Omega=\{1, x, x^*, xx^*, x^*x, xx^*x, x^*xx^*, \ldots\}.
\end{equation}
Moreover, denote by $\cP^*_k$ the space of $k$-homogeneous polynomial forms $P\in\cP^*$
and by $H^s_k$ the space of $\La^s(\bC^m)$-valued polynomial forms $P\in\cP^*_k$ which satisfy the Hodge-de Rham system (\ref{Hodge}).  
Then, using results from \cite{hom}, we shall deduce in Section 2 the corresponding Fischer decomposition for the $H$-action.

\begin{thm}\label{tFischer}
The space $\cP^*=\cP\otimes\La^*(\bC^m)$ decomposes as follows:
\begin{equation}\label{etFischer}
\cP^*=\cP^*_{(0,0)}\oplus\left(\bigoplus_{s=1}^{m-1}\bigoplus_{k=0}^{\infty}\cP^*_{(s,k)}\right)\oplus\cP^*_{(m,0)}\text{\ \ \ with\ \ \ }
\cP^*_{(s,k)}=\bigoplus_{w\in\Omega}wH^s_k.
\end{equation}
Moreover, in (\ref{etFischer}), all
$O(m)$-modules $H^s_k$ are non-trivial, irreducible and mutually inequivalent and all
$\cP^*_{(s,k)}$ are corresponding $O(m)$-isotypic components of $\cP^*.$
\end{thm}

Now we show that the Fischer decomposition of the space $\cP^*$ given in Theorem \ref{tFischer} is a~refinement of
the monogenic Fischer decomposition (\ref{FischerPin}).
Indeed, when we identify the Clifford algebra $\bC_m$ with the Grassmann algebra $\La^*(\bC^m)$ on the space $\cP^*$ we know that 
$$\underline{D}=d+d^*\text{\ \ \ and\ \ \ }-\underline{x}=x+x^*.$$
Consequently, the space of spherical monogenics of order $k$ is given by $$\cM_k=\{P\in\cP^*_k;\ (d+d^*)P=0\}.$$ 
Moreover, recall that the Laplace operator $\Delta$ is given by $\Delta=\sum_{j=1}^m\pa^2_{x_j}$ and put 
$$\Ker_k\Delta=\{P\in\cP^*_k;\ \Delta P=0\}.$$

As we mentioned before, under the $L$-action, $\bC_m\simeq\La^*(\bC^m)$ decomposes into many copies of basic spinor representations $\bS$ of $Pin(m)$ and so
the whole space $\cM_k$ of spherical monogenics  is not irreducible. Indeed,
under the $L$-action,
the space $\cM_k$
decomposes into many copies of irreducible modules $\cM_k(\bS)$.
In particular, we have thus that $\Ker_k\Delta=\cM_k\oplus(x+x^*)\cM_{k-1}$
and, by (\ref{Fischer}), we get easily the following decomposition of the space $\cP^*$
\begin{equation}\label{FischerPin*}
\cP^*=\bigoplus_{k=0}^{\infty}\bigoplus_{p=0}^{\infty}r^{2p}(\cM_k\oplus(x+x^*)\cM_{k-1}).
\end{equation}
In an obvious sense, the decomposition (\ref{FischerPin*}) is equivalent to (\ref{FischerPin}). 

In Section 3, we shall prove the following theorem which tells us that, under the $H$-action,
the spaces $\cM_k$ and $(x+x^*)\cM_{k-1}$ decompose again into many irreducible pieces 
but, in this case, these pieces are not equivalent and 
they have a~different representation character.

\begin{thm}\label{tdecompM}
For $k\geq 1,$ the following statements hold:
\begin{itemize}
\item[(a)] $\cM_k=\left(\bigoplus_{s=0}^m H^s_k\right)\oplus\left(\bigoplus_{s=1}^{m-1}M_{s,k}\right)$
\item[] where $M_{s,k}=[(k-1+m-s)x^*-(k-1+s)x]H^s_{k-1}.$ 
\item[(b)]
$(x+x^*)\cM_{k-1}=\left(\bigoplus_{s=0}^m(x+x^*)H^s_{k-1}\right)\oplus\left(\bigoplus_{s=1}^{m-1}W^s_k\right)$
\item[] where $W^s_k=[(k-2+m-s)xx^*-(k-2+s)x^*x]H^s_{k-2}.$
\end{itemize}
\end{thm}

Using Theorem \ref{tdecompM}, we get from the decomposition \eqref{FischerPin*} directly a~finer decomposition of the space $\cP^*$ which is irreducible with respect to the $H$-action.
 
The results stated in Theorems \ref{tFischer} and \ref{tdecompM} remain valid also for real valued polynomial forms, that is,
in the case when the complex Grassmann algebra $\La^*(\bC^m)$ is replaced with the real one $\La^*(\bR^m).$ 
Indeed, it is sufficient to realize that irreducible $O(m)$-representations $\La^s(\bC^m)$ are all of real type, see \cite[p. 163]{GM}.

 
\section{A proof of the Fischer decomposition for the $H$-action}\label{stFischer}

In this section, we give a~proof of Theorem \ref{tFischer} stated in Introduction. 
Let $\cP_k$ stand for the space of $k$-homogeneous (complex-valued) polynomials of $\cP$ and let
$\cP^s_k=\cP_k\otimes\Lambda^s(\bC^m).$
Then
it is easy to see that
\begin{equation}\label{decompP}
\cP^*=\bigoplus_{s=0}^m\bigoplus_{k=0}^{\infty}\cP^s_k.
\end{equation}
Let us remark that a~polynomial form $P$ belongs to $\cP^s_k$ if and only if
\begin{equation}
\label{skform}
P=\sum_I P_I\; dx_I
\end{equation}
where the sum is taken over all finite strictly increasing sequences $I=\{i_j\}_{j=1}^s$ of numbers of the set $\{1,\cdots, m\},$ $P_I\in\cP_k$ and $dx_I=dx_{i_1}\wedge\cdots\wedge dx_{i_s}.$
The contraction $dx_j\krat$ is defined as
$$dx_j\krat\; dx_I=\sum_{k=1}^s(-1)^{k-1}\delta_{ji_k}dx_{I\setminus\{i_k\}}\text{ and }
dx_j\krat\; P=\sum_I P_I\; dx_j\krat\; dx_I$$ for a~polynomial form $P.$ 
Then we have that
\begin{equation}
\label{dd^*}
d=\sum_{j=1}^m\; \pa_{x_j} dx_j\wedge\text{\ \ \ and\ \ \ }
d^*=-\sum_{j=1}^m\;\pa_{x_j} dx_j\krat, 
\end{equation}
\begin{equation}
\label{xx^*}
x=-\sum_{j=1}^m x_j\; dx_j\wedge\text{\ \ \ and\ \ \ }
x^*=\sum_{j=1}^m x_j\; dx_j\krat.
\end{equation}
It is easy to see that $d,$ $d^*,$ $x$ and $x^*$ are $O(m)$-invariant operators on the space $\cP^*.$

Now we describe explicitly an irreducible decomposition of $O(m)$-modules
$$\Ker_k^s\;\Delta=\{P\in\cP_k^s;\ \Delta P=0\}.$$
The following key result is obtained in \cite{hom}.

\begin{lem}{\label{thforms}}
Given $0\leq s\leq m$ and $k\in\NN_0,$ we have that
$$\Ker_k^s\Delta= H^s_k\oplus U^s_k\oplus V^s_k\oplus W^s_k$$ where
$H^s_k,$
$U^s_k,$ $V^s_k$ and $W^s_k$ are irreducible $O(m)$-modules
with the following properties:
\begin{itemize}
\item[(a1)] $H^s_k=\{P\in \cP^s_k;\ dP=0,\ d^*P=0\}$ and $\Ker_0^s\Delta= H^s_0=\cP_0^s.$
\item[(a2)] In addition, $H^s_k=\{0\}$ for $s\in\{0,m\}$ and $k\geq 1.$ Otherwise, all $O(m)$ modules $H^s_k$ are non-trivial, irreducible and mutually inequivalent.
\item[(b)]  $U^s_k=xH^{s-1}_{k-1}\simeq H^{s-1}_{k-1}$ for $1\leq s\leq m$ and $k\geq 1,$ and $U^s_k=\{0\}$ otherwise.
\item[(c)] $V^s_k=x^*H^{s+1}_{k-1}\simeq H^{s+1}_{k-1}$ for $0\leq s\leq m-1$ and $k\geq 1,$ and $V^s_k=\{0\}$ otherwise.
\item[(d)] $W^s_k=[(k-2+m-s)xx^*-(k-2+s)x^*x]H^s_{k-2}\simeq H^s_{k-2}$\\
for $1\leq s\leq m-1$ and $k\geq 2,$ and $W^s_k=\{0\}$ otherwise.
\end{itemize}
\end{lem}

Now we are ready to prove Theorem \ref{tFischer}. 

\begin{proof}[Proof of Theorem \ref{tFischer}]
As $\Ker^s_k\Delta=\cH_k\otimes\La^s(\bC^m)$ the Fischer decomposition (\ref{Fischer}) yields
$$\cP_k^s=\bigoplus_{p=0}^{[k/2]}r^{2p}\Ker^s_{k-2p}\Delta,$$
where for $y\in\bR,$ $[y]$ denotes the greatest integer not greater than $y.$
Consequently, by Lemma \ref{thforms}, we get the decomposition
\begin{equation}\label{decompPsk}
\cP_k^s= H^s_k\oplus
\bigoplus_{p=0}^{[k/2]}r^{2p}U^s_{k-2p}\oplus
\bigoplus_{p=0}^{[k/2]}r^{2p}V^s_{k-2p}\oplus
\bigoplus_{p=0}^{[k/2]}r^{2p}Z^s_{k-2p}
\end{equation}
where $Z^s_k=r^2 H^s_{k-2}\oplus W^s_k.$
Since $r^2=-(xx^*+x^*x)$ Lemma \ref{thforms} implies that,
for $0\leq s\leq m$ and $k\geq 2,$   
$$Z^s_k=(xx^*)H^s_{k-2}\oplus (x^*x)H^s_{k-2}.$$
Moreover, it is easy to see that 
$$r^{2p}U^s_k=(xx^*)^pxH^{s-1}_{k-1},\ \ \
r^{2p}V^s_k=(x^*x)^px^*H^{s+1}_{k-1}\text{\ \ and}$$
$$r^{2p}Z^s_k=(xx^*)^{p+1}H^s_{k-2}\oplus (x^*x)^{p+1}H^s_{k-2}.$$
Now to complete the proof
it suffices to use the decompositions (\ref{decompP}) and (\ref{decompPsk}).
\end{proof}

At the end of this section we collect the well-known relations we need later on. 
Put, for linear operators $T_1$ and $T_2$ on $\cP^*,$
$\{T_1,T_2\}=T_1T_2+T_2T_1$ and $[T_1,T_2]=T_1T_2-T_2T_1.$  Then we have that (see e.g. \cite{hom} or \cite{BDS})

\begin{lem}\label{lrels} Let $E$ be the Euler operator and $\hat{E}$ be the skew Euler operator, i.e. $$E=\sum_{j=1}^m x_j\pa_{x_j}\text{\ \ \ and\ \ \ }
\hat{E}=\sum_{j=1}^m(dx_j\;\wedge)(dx_j\krat).$$
Then we have that $EP=kP$ and $\hat{E}P=sP$ for each $P\in\cP^s_k.$

Furthermore, the following relations hold:
\begin{equation*}
\begin{array}{lll}
\{x,x\}=0, &\{x^*,x^*\}=0, &\{x,x^*\}=-r^2,\medskip\\{}
\{d,d\}=0, &\{d^*,d^*\}=0, &\{d,d^*\}=-\Delta,\medskip\\{}
\{x^*,d\}=E+\hat{E}, &\{x,d^*\}=E-\hat{E}+m, &\{x^*,d^*\}=0=\{x,d\}.
\end{array}
\end{equation*}
\end{lem}

Using Lemma \ref{lrels}, we may give, for example, an explicit description of the projections of the space $\Ker_k^s\Delta$ onto the pieces $H^s_k,$ $U^s_k,$ $V^s_k$ and $W^s_k.$ 

\begin{prop}\label{chforms}
Given $0\leq s\leq m$ and $k\in\NN_0,$ put $c_1=k-2+s$ and $c_2=k-2+m-s.$
Furthermore, let $\pi_1,$ $\pi_2,$ $\pi_3$ and $\pi_4$ be the projections of the space $\Ker_k^s\Delta$ onto the subspaces $H^s_k,$
$U^s_k,$ $V^s_k$ and $W^s_k,$ respectively. Then we have that
\begin{equation*}
\pi_4=\left\{
\begin{array}{ll}
\frac{c_2xx^*-c_1x^*x}{c_1c_2(c_1+c_2+2)}\;dd^* &\text{for\ \ }1\leq s\leq m-1\text{\ and\ }k\geq 2,\medskip\\{}
0, &\text{otherwise.}
\end{array}
\right.
\end{equation*}
Moreover, denoting $\pi=1-\pi_4,$ we have that
\begin{equation*}
\begin{array}{ll}
&\pi_2 =\left\{
\begin{array}{ll}
\frac{1}{c_2+2}\;xd^*\pi &\text{for\ \ }1\leq s\leq m\ \text{\ and\ }k\geq 1,\medskip\\{}
0, &\text{otherwise;}
\end{array}
\right.\medskip\\{}

&\pi_3 =\left\{
\begin{array}{ll}
\frac{1}{c_1+2}\;x^*d\pi &\text{for\ \ }0\leq s\leq m-1\text{\ and\ }k\geq 1,\medskip\\{}
0, &\text{otherwise;}
\end{array}
\right.\medskip\\{}

&\pi_1 =1-\pi_2-\pi_3-\pi_4.
\end{array}
\end{equation*} 
\end{prop}

\begin{proof}
Let $P\in\Ker^s_k\Delta$ be given. Then, by Lemma \ref{thforms}, there are uniquely determined $P_1\in H^s_k,$ $P_2\in H^{s-1}_{k-1},$ $P_3\in H^{s+1}_{k-1}$ and $P_4\in H^s_{k-2}$ such that $$P=P_1+xP_2+x^*P_3+(c_2xx^*-c_1x^*x)P_4.$$
By Lemma \ref{lrels},
it is easy to see that
$$dd^*P=c_1c_2(c_1+c_2+2)P_4=-d^*dP,$$ which easily implies the formula for $\pi_4.$

Moreover, $\pi(P)=P_1+xP_2+x^*P_3.$ By Lemma \ref{lrels}, we have that 
$$d^*\pi(P)=d^*xP_2=(c_2+2)P_2\text{\ \ \ and\ \ \ }d\pi(P)=dx^*P_3=(c_1+2)P_3,$$
from which the formulae for the projections $\pi_2$ and $\pi_3$ may be derived.
\end{proof}

\section{Decomposition of monogenic polynomial forms}
 
In this section, we give a~proof of Theorem \ref{tdecompM} stated in Introduction.
To prove Theorem \ref{tdecompM} we need some lemmas.

\begin{lem}\label{lMsk}
For $1\leq s\leq m-1$ and $k\geq 1,$
we have that
$$\left(xH^s_{k-1}\oplus x^*H^s_{k-1}\right)\cap\cM_k=M_{s,k}.$$
Here $M_{s,k}=[(k-1+m-s)x^*-(k-1+s)x]H^s_{k-1}.$
\end{lem}

\begin{proof}
Let $P_1,P_2\in H^s_{k-1}$ and put $P=xP_1+x^*P_2.$ It suffices to show that $(d+d^*)P=0$
if and only if $$P_1=-\frac{k-1+s}{k-1+m-s}P_2.$$
By virtue of Lemma \ref{lrels}, it is easy to see that
$$(d+d^*)P=(E+m-\hat E)P_1+(E+\hat E)P_2=(k-1+m-s)P_1+(k-1+s)P_2,$$
which completes the proof.
\end{proof}

\begin{lem}\label{ldUV}
For $1\leq s\leq m-1$ and $k\geq 1,$ we have that
$$xH^s_{k-1}\oplus x^*H^s_{k-1}=(x+x^*)H^s_{k-1}\oplus M_{s,k}. $$
\end{lem}

\begin{proof}
Obvious.
\end{proof}

\begin{proof}[Proof of Theorem \ref{tdecompM}]
Put$$\tilde\cM_k=\left(\bigoplus_{s=0}^m H^s_k\right)\oplus\left(\bigoplus_{s=1}^{m-1}M_{s,k}\right).$$
Then, by Lemma \ref{lMsk}, it is easy to see that
$\tilde\cM_k\subset\cM_k.$
Moreover, by Lemma \ref{thforms}, $W^s_k=(x+x^*)M_{s,k-1}.$
Finally, using Lemma \ref{ldUV} and Lemma \ref{thforms},
we obtain that
$$
\Ker_k\Delta=\tilde\cM_k\oplus(x+x^*)\tilde\cM_{k-1}\subset\cM_k\oplus(x+x^*)\cM_{k-1}=\Ker_k\Delta,
$$
which completes the proof.
\end{proof}


\subsection*{Acknowledgments}

R. L\'avi\v cka and V. Sou\v cek acknowledge the financial support from the grant GA 201/08/0397.
This work is also a part of the research plan MSM 0021620839, which is financed by the Ministry of Education of the Czech Republic.



\bigskip\bigskip\bigskip

\noindent
Richard Delanghe,\\ Clifford Research Group, Department of Mathematical Analysis,\\ Ghent University, Galglaan 2, B-9000 Gent, Belgium\\
email: \texttt{richard.delanghe@ugent.be}

\bigskip

\noindent
Roman L\'avi\v cka and Vladim\'ir Sou\v cek,\\ Mathematical Institute, Charles University,\\ Sokolovsk\'a 83, 186 75 Praha 8, Czech Republic\\
email: \texttt{lavicka@karlin.mff.cuni.cz} and \texttt{soucek@karlin.mff.cuni.cz}


\begin{thebibliography}{99}
\def\topsep{0pt}
\def\parsep{0pt plus 5pt minus 1pt}
\def\itemsep{-0.5ex} 
\small
\bibitem{BGLS} S. Bock, K. G\"urlebeck, R. L\' avi\v cka and V. Sou\v cek, The Gelfand-Tsetlin bases for spherical monogenics in dimension 3, preprint.
%
 \bibitem{Bock2010c} S. Bock, Orthogonal Appell bases in dimension 2,3 and 4. In Numerical Analysis and Applied Mathematics
(T.E. Simos, G. Psihoyios, and Ch. Tsitouras, eds.), AIP Conference Proceedings, vol. 1281. American
Institute of Physics: Melville, NY, 2010; 1447--1450.
%
\bibitem{Bock2010a} S. Bock, On a three dimensional analogue to the holomorphic $z$-powers: Power series and recurrence
formulae, submitted, 2010.
%
\bibitem{Bock2009} S. Bock, \"{U}ber funktionentheoretische {M}ethoden in der r\"{a}umlichen
{E}lastizit\"{a}tstheorie, PhD thesis, Bauhaus-University, Weimar, (url:
http://e-pub.uni-weimar.de/frontdoor.php?source\_opus=1503, date: 07.04.2010), 2009.
%
\bibitem{BG} S. Bock and K. G\"urlebeck, On a~generalized Appell
system and monogenic power series, Mathematical Methods in the Applied Sciences 33 (2010),
394--411.
%
\bibitem{BDS} F. Brackx, R. Delanghe and F. Sommen, Differential forms and/or multi-vector
functions, CUBO 7 (2005), 139-170.
%
\bibitem{BSES} F. Brackx, H. De Schepper, D. Eelbode and V. Sou\v cek, The Howe dual pair in hermitian Clifford analysis, Rev. Mat. Iberoamericana 26 (2010)(2), 449-479.
%
\bibitem{ckH} F.\ Brackx, H.\ De Schepper, R.\ L\'{a}vi\v{c}ka, V.\ Sou\v{c}ek,
\textit{The Cauchy-Kovalevskaya	Extension Theorem in Hermitean Clifford Analysis}, preprint.
%
\bibitem{kerH} F.\ Brackx, H.\ De Schepper, R.\ L\'{a}vi\v{c}ka, V.\ Sou\v{c}ek, 
\textit{Fischer decompositions of kernels of Hermitean Dirac operators}, In: T.E.\ Simos, G.\ Psihoyios, Ch.\ Tsitouras, {\em Numerical Analysis and Applied Mathematics}, AIP Conference Proceedings, Rhodes, Greece (2010).
%
\bibitem{GTinH} F.\ Brackx, H.\ De Schepper, R.\ L\'{a}vi\v{c}ka, V.\ Sou\v{c}ek, 
\textit{Gel'fand-Tsetlin procedure for the construction of orthogonal bases in Hermitean Clifford analysis},
In: T.E.\ Simos, G.\ Psihoyios, Ch.\ Tsitouras, {\em Numerical Analysis and Applied Mathematics}, AIP Conference Proceedings, Rhodes, Greece (2010).
%
\bibitem{GT2H} F.\ Brackx, H.\ De Schepper, R.\ L\'{a}vi\v{c}ka, V.\ Sou\v{c}ek, 
\textit{Orthogonal basis of Hermitean monogenic polynomials: an explicit construction in complex dimension $2$}, In: T.E.\ Simos, G.\ Psihoyios, Ch.\ Tsitouras (eds.), {\em Numerical Analysis and Applied Mathematics}, AIP Conference Proceedings, Rhodes, Greece (2010). 
%
 \bibitem{CacGueMal} I. Ca\c c\~ao, K. G\"urlebeck, H.R. Malonek, Special monogenic polynomials and $L_2$-approximation. Advances in Applied Clifford Algebras 2001; 11(S2):47--60.
%
\bibitem{CacGueBock} I. Ca\c c\~ao, K. G\"urlebeck, S. Bock. Complete orthonormal systems of spherical monogenics - a constructive approach. In Methods of Complex and
Clifford Analysis, Son LH, Tutschke W, Jain S (eds). Proceedings of ICAM, Hanoi, SAS International Publications, 2004.
%
 \bibitem{BockCacGue} I. Ca\c c\~ao, K. G\"urlebeck, S. Bock, On derivatives of spherical monogenics. Complex Variables and Elliptic Equations 2006; 51(811):847--869.
%
 \bibitem{CM06} I. Ca\c c\~ao and H.~R. Malonek, Remarks on some properties of monogenic polynomials,
ICNAAM 2006. International conference on numerical analysis and applied mathematics 2006 (T.E.
Simos, G. Psihoyios, and Ch. Tsitouras, eds.), Wiley-VCH, Weinheim, 2006, pp. 596-599.
%
 \bibitem{CM08} I. Ca\c c\~ao and H.~R. Malonek, On a complete set of hypercomplex Appell polynomials,
Proc. ICNAAM 2008, (T. E. Timos, G. Psihoyios, Ch. Tsitouras, Eds.), AIP Conference Proceedings
1048, 647-650.
%
 \bibitem{cac} I. Ca\c c\~ao, \textit{Constructive approximation by monogenic polynomials}, PhD thesis, Univ. Aveiro, 2004.
%
\bibitem{DLS} R. Delanghe, R. L\'avi\v cka and V. Sou\v cek, On polynomial solutions of generalized Moisil-Th\'eodoresco systems and Hodge-de Rham systems, to appear in Adv.
appl. Clifford alg. (arXiv:0908.0842 [math.CV], 2009).
%
\bibitem{DLS3} R. Delanghe, R. L\'avi\v cka and V. Sou\v cek, The Howe duality for Hodge systems, In: Proceedings of 18th International Conference on the Application of
Computer Science and Mathematics in Architecture and Civil Engineering
(ed. K. G\"urlebeck and C. K\"onke), Bauhaus-Universit\"at Weimar, Weimar,
2009.
%
\bibitem{DLS4} R. Delanghe, R. L\'avi\v cka and V. Sou\v cek, The Gelfand-Tsetlin bases for Hodge-de Rham systems in Euclidean spaces, preprint.
%
\bibitem{FCM} M.~I. Falc\~ao, J.~F. Cruz and H.~R. Malonek, Remarks on the generation of monogenic
functions, Proc. of the 17-th International Conference on the Application of Computer Science and
Mathematics in Architecture and Civil Engineering, ISSN 1611-4086 (K. G\"urlebeck and C. K\"onke,
eds.), Bauhaus-University Weimar, 2006.
%
\bibitem{FM} M.~I. Falc\~ao and H.~R. Malonek, Generalized exponentials through Appell sets in $\bR^{n+1}$ and
Bessel functions, Numerical Analysis and Applied Mathematics (T.E. Simos, G. Psihoyios, and Ch.
Tsitouras, eds.), AIP Conference Proceedings, vol. 936, American Institute of Physics, 2007, pp.
750-753 (ISBN: 978-0-7354-0447-2).

\bibitem{FSS} L. Frappat, P. Sorba and A. Sciarrino, Dictionary on Lie superalgebras, arXiv:hep-th/960761v1, 1996.
%
\bibitem{GM} J. E. Gilbert and M. A. M. Murray, Clifford Algebras and Dirac Operators in Harmonic Analysis,
Cambridge University Press, Cambridge, 1991.
%
\bibitem{Goo} R. Goodman, Multiplicity-free spaces and Schur-Weyl-Howe duality, In: Representations of real and $p$-adic groups (E. C. Tan and C. B. Zhu eds.), Lecture note series - Institute for mathematical sciences, Vol. 2, World scientific, Singapore, 2004.
%
\bibitem{GW} R. Goodman and N. Wallach, Representations and invariants of the classical groups, Cambridge University Press, Cambridge, 1998.
%
\bibitem{gs1} K. G\"urlebeck, W. Spr\"ossig, {\it Quaternionic analysis and boundary value problems}, Birkh\"auser, Basel, 1990.
%
\bibitem{gs2} K. G\"urlebeck, W. Spr\"ossig, {\it Quaternionic and Clifford calculus for physicists and engineers}, Wiley,
Chichester, 1997.
%
\bibitem{Gurlebeck1999} K. G\"{u}rlebeck and H.~R. Malonek, A hypercomplex derivative of monogenic functions
in $\mathbb{R}^{n+1}$ and its applications, Complex Variables 39 (1999), 199--228.
%
\bibitem{NGue2009} N. G\"urlebeck, On Appell Sets and the Fueter-Sce Mapping, Advances in Applied Clifford Algebras 19 (2009), 51-61.
%
\bibitem{hom} Y. Homma, Spinor-valued and Clifford algebra-valued harmonic polynomials, J. Geom. Phys. 37 (2001), 201-215.
%
\bibitem{how} R. Howe, Remarks on classical invariant theory, Trans. Am. Math. Soc. 313 (1989) (2), 539-570.
%
\bibitem{lav_Fischer} R. L\'avi\v cka, The Fischer Decomposition for the H-action and Its Applications, arXiv:1002.0527v1 [math.CV], 2010, to appear. 
%
\bibitem{lavSL2} R.~L\'avi\v cka, Canonical bases for sl(2,C)-modules of spherical monogenics in dimension 3, arXiv:1003.5587v2 [math.CV], 2010,
to appear in Arch. Math.
%
\bibitem{lavSL3} R.~L\'avi\v cka, 
On the Structure of Monogenic Multi-Vector Valued Polynomials, In: ICNAAM 2009, Rethymno, Crete, Greece, 18-22 September 2009 (eds. T. E. simos, G. Psihoyios and Ch. Tsitouras), AIP Conf. Proc. 1168 (2009)(793), pp. 793-796.
 %
 \bibitem{step2} R.~L\'avi\v cka, V.~Sou\v cek, P.~Van Lancker, Spherical monogenics: step two branching, preprint.
%
\bibitem{CM07} H.~R. Malonek, M.~I. Falc\~ao, Special monogenic polynomials\,-\,properties and applications. In Numerical Analysis and Applied Mathematics (T.E.
Simos, G. Psihoyios, and Ch. Tsitouras, eds.), AIP Conference Proceedings, vol. 936. American Institute of
Physics: Melville, NY, 2007; 764--767.
%
\bibitem{Malonek1987} H.~R. Malonek, {Z}um {H}olomorphiebegriff in h\"{o}heren {D}imensionen, Habilitationsschrift. P\"{a}dagogische Hochschule Halle,
 1987.
%
\bibitem{MPS} H.R. Malonek, D. Pe\~na Pe\~na and F. Sommen, Fischer decomposition by inframonogenic
functions, arXiv:0911.0070 [math.CV], 2009 (to appear in CUBO).
%
\bibitem{MR} H.R. Malonek and G. Ren, Almansi-type theorems in Clifford analysis, Math. Meth. Appl. Sci. 25 (2002), 1541-1552.
%
 \bibitem{MS} I.~M. Mitelman and M.~V. Shapiro, Differentiation of
the Martinelli--Bochner integrals and the notion of hyperderivability. {\em Math. Nachr.} 172:
(1995), 211--238.
%
\bibitem{mor09} J. Morais, Approximation by homogeneous polynomial solutions of the Riesz system in $\bR^3$, PhD thesis, Bauhaus-Univ., Weimar, 2009.
%
\bibitem{rya} J. Ryan, Iterated Dirac operators in $\bC^n,$ Z. Anal. Anwendungen 9 (1990), 385-401.
%
\bibitem{som} F. Sommen, Spingroups and spherical means III, Rend. Circ. Mat. Palermo
(2) Suppl. No 1 (1989), 295-323.
%
\bibitem{S_Maxwell} W. Spr\"ossig, Quaterionic analysis and Maxwell's equations, CUBO A mathematical journal 7 (2005) (2), 57-67.
%
\bibitem{spr} W. Spr\"ossig, On Helmholtz decompositions and their generalizations - An overview, Math. Meth. Appl. Sci. 33 (2010), 374-383.
%
\bibitem{van} P. Van Lancker, Spherical Monogenics: An Algebraic Approach, Adv. appl. Clifford alg. 19 (2009), 467-496.
%
\end{thebibliography}
\end{document}